\documentclass[10pt]{article}
\usepackage{latexsym}
\usepackage{amsfonts}
\usepackage{enumerate}
\usepackage{multicol}
\usepackage{graphicx}
\usepackage{amssymb}
\usepackage{amsmath}
\usepackage{epic}

\title{On Perfect Bases in Finite Abelian Groups \\[.4in]}

\author{B\'{e}la Bajnok\footnote{Department of Mathematics, Gettysburg College, U.S.A.  Email: bbajnok@gettysburg.edu. Corresponding author.} 
\hspace{.7in} Connor Berson\footnote{Asymmetrik, U.S.A.  Email: cjberson@gmail.com.} 
\hspace{.7in} Hoang Anh Just\footnote{Department of Electrical and Computer Engineering, Virginia Polytechnic Institute and State University, U.S.A. E-mail:  just@vt.edu.}}

\date{September 29, 2021}

\newtheorem{thm}{Theorem}

\newtheorem{prop}[thm]{Proposition}

\newtheorem{que}[thm]{Question}

\begin{document}

\maketitle

\begin{abstract}

Let $G$ be a finite abelian group and $s$ be a positive integer.  A subset $A$ of $G$ is called a {\em perfect $s$-basis of $G$} if each element of $G$ can be written uniquely as the sum of at most $s$ (not-necessarily-distinct) elements of $A$; similarly, we say that $A$ is a {\em perfect restricted $s$-basis of $G$} if each element of $G$ can be written uniquely as the sum of at most $s$ distinct elements of $A$.  We prove that perfect $s$-bases exist only in the trivial cases of $s=1$ or $|A|=1$.  The situation is different with restricted addition where perfection is more frequent; here we treat the case of $s=2$ and prove that $G$ has a perfect restricted $2$-basis if, and only if, it is isomorphic to $\mathbb{Z}_2$, $\mathbb{Z}_4$, $\mathbb{Z}_7$, $\mathbb{Z}_2^2$, $\mathbb{Z}_2^4$, or $\mathbb{Z}_2^2 \times \mathbb{Z}_4$.


      
\end{abstract}

\thispagestyle{empty}

\section{Introduction}

Throughout this paper we let $G$ be an additively written finite abelian group.  When $G$ is cyclic and of order $n$, we identify it with $\mathbb{Z}_n=\mathbb{Z}/n\mathbb{Z}$; we consider $0,1,\dots,n-1$ interchangeably as integers and as elements of $\mathbb{Z}_n$.

For subsets $A$ and $B$ of $G$, we let $A+B$ and $A-B$ denote the set of two-term sums and differences, respectively, with one term chosen from $A$ and the other from $B$.  
(If, say, $A$ consists of a single element $a$, we simply write $a+B$ and $a-B$ instead of $A+B$ and $A-B$.)  
For a subset $A$ of $G$ and a nonnegative integer $h$, $hA$ denotes the {\em $h$-fold  sumset} of $A$, that is, the collection of $h$-term sums with (not-necessarily-distinct) elements from $A$.  The study of sumsets is the main theme in much of additive number theory;  Nathanson's textbook \cite{Nat:1996a} provides a general introduction to the field, while the monograph \cite{Baj:2018a} discusses and provides the background for over 300 of its open questions.      

Given a positive integer $s$, we say that a subset $A$ of $G$ is an {\em $s$-basis} for $G$ if every element of $G$ can be written as the sum of at most $s$ (not-necessarily-distinct) elements of $A$; that is, if $\cup_{h=0}^s hA=G$.  The concept has enjoyed a rich history since it was first discussed by Erd\H{o}s and Tur\'an  in \cite{ErdTur:1941a} eighty years ago; see, for example, \cite{ErdGra:1980b}, \cite{ErdNat:1987a}, \cite{LamThaPla:2017a}, \cite{Nat:1974a},  \cite{Nat:1996a}, and \cite{Nat:2014a}.

A subset $A$ of $G$ is called a {\em perfect $s$-basis} if every element of $G$ can be written uniquely (apart from the order of the terms) as a sum of at most $s$ elements of $A$.  (In somewhat of a contrast to our usual understanding of bases, this uniqueness property is not part of the definition of $s$-bases where only the spanning property is required.) One can then ask for each instance when a group $G$ possesses a perfect $s$-basis.  
Trivially, the set of nonzero elements is a perfect $1$-basis in any group $G$, and the 1-element set consisting of a generator of $G$ is a perfect $s$-basis in the cyclic group of order $s+1$. 
We prove that there are no others:

\begin{thm} \label{CoHo-1}
If a subset $A$ of a finite abelian group $G$ is a perfect $s$-basis, then $s=1$ and $A=G \setminus \{0\}$, or $G\cong \mathbb{Z}_{s+1}$ and $|A|=1$.

\end{thm}
Our proof of Theorem \ref{CoHo-1} is based on the fact that if $A$ is a perfect $s$-basis in $G$, then the $s+1$ subsets 
$$-A, \; A-A, \; 2A-A, \; \dots, \; (s-1)A-A, \; \mbox{and} \; (s-1)A$$ are pairwise disjoint.  This will allow us to conclude that, when $A$  is a perfect $s$-basis in $G$ with $|A| \geq 2$ and $s \geq 2$, then the set
$$S=\cup_{h=0}^{s-1} (hA-A) \; \cup \; (s-1)A$$ has size more than the order of $G$, which is clearly not possible.    

We should note that focusing on sums of at most $s$ terms is equivalent to the case when exactly $s$ terms are considered.  A subset $B$ of $G$ is said to be an {\em additive basis of order $s$} if every element of $G$ can be expressed as the sum of exactly $s$ terms (that is, if $sA=G$), and it is called a {\em $B_s$ set} if all $s$-term sums are distinct (up to the order of the terms).  We now prove that $G$ has a perfect $s$-basis of size $m$ if, and only if, it has a subset of size $m+1$ that is an additive basis of order $s$ as well as a $B_s$ set.

Suppose first that $A$ is a perfect $s$-basis in $G$.  Since then $0 \not \in A$, the set $B=A \cup \{0\}$ has size $|A|+1$; it is also easy to see that each element of $G$ can be written uniquely as the sum of exactly $s$ elements of $B$.  Conversely, assume that $B$ is an additive basis of order $s$ as well as a $B_s$ set.  Note that this implies that each translate of $B$ has these two properties as well; in particular, we may assume that $0 \in B$.  This then clearly implies that $A=B \setminus \{0\}$ is a perfect $s$-basis in $G$.  We can thus restate Theorem \ref{CoHo-1} as follows:

\begin{thm} \label{CoHo-2}
If a subset $B$ of a finite abelian group $G$ is a $B_s$ set and an additive basis of order $s$, then $s=1$ and $B=G$, or $G\cong \mathbb{Z}_{s+1}$ and $|B|=2$.

\end{thm}

Let us now turn to {\em restricted addition}, that is, when the terms in the sums must be pairwise distinct.  The {\em $h$-fold restricted sumset} of $A$, denoted by $h \hat{\;}A$, is the collection of $h$-term sums with pairwise distinct elements from $A$.  
  We call a subset $A$ of $G$ a {restricted \em $s$-basis} if every element of $G$ can be written as the sum of at most $s$ distinct elements of $A$, and $A$ is said to be a {\em perfect restricted $s$-basis} if these sums are (apart from the order of the terms) unique  for each element of the group.  The instances of $|A|=1$ or $s=1$ are identical for restricted addition and unrestricted addition, and are as listed above.  As a major contrast, however, there are infinitely many perfect restricted $s$-bases with $|A| \geq 2$ and $s \geq 2$ (see Section 4).

Our strategy for classifying perfect restricted $s$-bases is similar to the unrestricted case, but exhibiting sets that are pairwise disjoint seems more elusive.  Here we only treat the case of $s=2$, for which we prove that when $G$ is not isomorphic to the elementary abelian 2-group, then $G$ has a perfect restricted $2$-basis in exactly three cases: when $G$ is isomorphic to $\mathbb{Z}_4$, $\mathbb{Z}_7$, or $\mathbb{Z}_2^2 \times \mathbb{Z}_4$.  Similarly to the unrestricted case, for a given perfect restricted $2$-basis $A$ we consider the set
$$T=(A-A) \cup A,$$ and prove that, unless $G$ is isomorphic to the three groups just listed, $T$ has size more than the order of $G$.      
Elements of order $2$ in $G$ require additional attention that we are able to handle when $G$ is not isomorphic to an elementary abelian 2-group.  However, if all nonzero elements of $G$ have order 2, then $$A-A=2A=\{0\} \cup 2 \hat{\;} A,$$ and thus $T=\cup_{h=0}^2 h\hat{\;} A.$ Therefore, $A$ being a perfect restricted $2$-basis in $G$ is equivalent to having $T=G$, resulting in no contradiction.  Luckily, a problem of Ramanujan comes to the rescue.

In 1913, Ramanujan asked in \cite{Ram:2000a} whether the quantity $2^k-7$ can be a square number for any integer $k$ besides 3, 4, 5, 7, and 15 (see also Question 464 in \cite{BerChoKan:1999a}).   The negative answer was given by Nagell in 1948 (see \cite{Nag:1948a}; \cite{Nag:1961a} for the English version).  Suppose now that $G$ is the elementary abelian 2-group of rank $r$ that has a perfect restricted $2$-basis of size $m$: we then must have 
$$2^r={m \choose 0} + {m \choose 1} + {m \choose 2},$$ and therefore $$2^{r+3}-7=(2m+1)^2.$$  We thus have exactly four choices for $r$: 1, 2, 4, and $12$.  We show that perfect restricted $2$-bases exist in the first three cases but not in $\mathbb{Z}_2^{12}$.   In summary, we have the following results:  
  
\begin{thm} \label{CoHo-3}
A finite abelian group $G$ has a perfect restricted $2$-basis if, and only if, it is isomorphic to one of the following groups: $\mathbb{Z}_2, \; \mathbb{Z}_4, \; \mathbb{Z}_7, \; \mathbb{Z}_2^2, \; \mathbb{Z}_2^4$, or $\mathbb{Z}_2^2 \times \mathbb{Z}_4$.

\end{thm}

The structure of the paper is as follows: In Section 2 we prove Theorem \ref{CoHo-1}, and in Section 3 we prove Theorem \ref{CoHo-3}.  We then consider some further results and present some open questions in Section 4.

\section{The proof of Theorem \ref{CoHo-1}}

As mentioned in the Introduction, in order to prove Theorem  \ref{CoHo-1}, we argue indirectly, and show that, when $A$  is a perfect $s$-basis in $G$ with $|A| \geq 2$ and $s \geq 2$, then  the size of $$S=\cup_{h=0}^{s-1} (hA-A) \; \cup \; (s-1)A$$ is more than the order of $G$, which  is impossible.

Recall that if $A$ is a perfect $s$-basis in $G$, then it is a $B_h$ set for every $0 \leq h \leq s$, and $G$ is the disjoint union of sets $hA$ with $0 \leq h \leq s$.  Therefore, for each such $h$ we have $$|hA|={m+h-1 \choose h}$$ and thus
$$|G|=\sum_{h=0}^s |hA|= {m+s \choose s}.$$ 

In order to compute the size of $S$, first note that if $A$ is a perfect $s$-basis in $G$, then the $s+1$ sets $-A, A-A, 2A-A, \dots, (s-1)A-A$, and $(s-1)A$ are pairwise disjoint.  Therefore, 
$$|S|=\sum_{h=0}^{s-1} |hA-A| + |(s-1)A|= \sum_{h=0}^{s-1} |hA-A| + {m+s-2 \choose s-1}.$$

For $h=0$ we have $|hA-A|=|-A|=m$. To find the size of $hA-A$ for a given $1 \leq h \leq s-1$, we first observe that $(h-1)A$ is a subset of $hA-A$, so there are exactly ${m+h-2 \choose h-1}$ elements of $hA-A$ that are in $(h-1)A$.  We will now show that there are precisely $m \cdot {m+h-2 \choose h}$ elements in $hA-A$ that are not in $(h-1)A$.

Let $A=\{a_1, \dots, a_m \}$.  Observe that each element $g$ of $(hA-A) \setminus (h-1)A$ can be written as 
$$g=\lambda_1 a_1+ \cdots + \lambda_m a_m - a_i,$$ where $\lambda_1, \dots, \lambda_m$ are nonnegative integers that sum to $h$, $1 \leq i \leq m$, and $\lambda_i=0$.  As we now show, the expression is unique for each $g$.  Suppose that we also have   
$$g=\lambda'_1 a_1+ \cdots + \lambda'_m a_m - a_j$$  with some $1 \leq j \leq m$, $\lambda'_1+ \cdots + \lambda'_m=h$, and $\lambda'_j=0$.  This yields  
$$\lambda_1 a_1+ \cdots + \lambda_m a_m + a_j = \lambda'_1 a_1+ \cdots + \lambda'_m a_m + a_i \in (h+1)A.$$  If $i \neq j$, then   
 since $A$ is a $B_{h+1}$ set, we must have $\lambda_j+1 = \lambda'_j$ and $\lambda_i = \lambda'_i+1$, contradicting $\lambda_i = \lambda'_j=0$.
Furthermore, $g+a_i \in hA$, so the coefficients  $\lambda_1, \dots, \lambda_m$ are unique as well.

This means that 
$(hA-A) \setminus (h-1)A$ is the disjoint union of sets $h(A \setminus \{a_i\}) -a_i$ for $1 \leq i \leq m$.  Since a subset of a $B_h$ set is also a $B_h$ set, we have 
$$|h(A \setminus \{a_i \}) -a_i|=|h(A \setminus \{a_i\})|={m+h-2 \choose h},$$
implying $$|(hA-A) \setminus (h-1)A |= m \cdot {m+h-2 \choose h},$$ as claimed.

We thus find that $$|hA-A|=m \cdot {m+h-2 \choose h} + {m+h-2 \choose h-1},$$ from which
$$\sum_{h=0}^{s-1} |hA-A| = m+ \sum_{h=1}^{s-1} m \cdot {m+h-2 \choose h} + \sum_{h=1}^{s-1}{m+h-2 \choose h-1} = m {m+s-2 \choose s-1} + {m+s-2 \choose s-2}.$$  This then yields
$$|S|=\sum_{h=0}^{s-1} |hA-A| + |(s-1)A|= m {m+s-2 \choose s-1} + {m+s-1 \choose s-1};$$ 
rewriting this expression, we get
$$|S|= m {m+s-2 \choose s-1}  + {m+s \choose s} - {m+s-1 \choose s} = {m+s \choose s} + \frac{(m-1)(s-1)}{s} {m+s-2 \choose s-1}.$$
Since this quantity is larger than $|G|={m+s \choose s}$ for $m \geq 2$ and $s \geq 2$, our proof is complete.

\section{The proof of Theorem \ref{CoHo-3}}

We start by treating the cases when $G$ is isomorphic to an elementary abelian 2-group $\mathbb{Z}_2^r$.  As explained in the Introduction, we only need to examine four cases: when $G$ has rank 1, 2, 4, or 12.  It is not hard to find perfect restricted 2-bases in the first three of these cases; for example:

\begin{itemize}
 \item in $\mathbb{Z}_2$, the set $\{1\}$ is a perfect restricted 2-basis;

 \item in $\mathbb{Z}_2^2$, the set  $\{(0,1),(1,0)\}$ is a perfect restricted 2-basis; and

 \item in $\mathbb{Z}_2^4$, the set  $\{(0, 0, 0, 1), (0, 0, 1, 0), (0, 1, 0, 0), (1, 0, 0, 0), (1, 1, 1, 1)\}$ is a perfect restricted 2-basis.

\end{itemize}  However, we can easily see that $\mathbb{Z}_2^{12}$ has no perfect restricted 2-bases.  For the sake of a contradiction, suppose that $A$ is a perfect restricted 2-basis in $\mathbb{Z}_2^{12}$; it then has size $|A|=90$.  If $k$ elements of $A$ have a first component of $1$ for some $0 \leq k \leq 90,$ then $0\hat{\;}A \cup 1\hat{\;}A \cup 2 \hat{\;}A$ has $k+k(90-k)$ elements with a first component of $1$.  Since there is no $k$ for which $k+k(90-k)$ equals $|G|/2=2^{11}$, this is a contradiction.

We now turn to the cases when $G$ is not isomorphic to elementary abelian 2-groups: we need to prove that $G$ has a perfect restricted 2-basis if, and only if, it is isomorphic to $\mathbb{Z}_4$, $\mathbb{Z}_7$, or $\mathbb{Z}_2^2 \times \mathbb{Z}_4$.  

We let $L$ denote the subset consisting of the identity element of $G$ as well as all involutions in $G$; that is, $$L=\{g \in G \mid 2g=0\}.$$
Note that $L$ is a subgroup of $G$; in fact, $L$ is
isomorphic to the elementary abelian 2-group whose rank equals the number of even-order
terms in the invariant decomposition of $G$.
 
Let $A=\{a_1,\dots,a_m\}$ be an $m$-subset of $G$, and suppose that $A$ is a perfect restricted 2-basis in $G$.  The main idea of the proof is to estimate the size of $(A-A) \cup A$, namely, to prove that
$$|(A-A) \cup A| \geq m^2-2m+2-|L|.$$ Since, as we show, this bound is more than $|G|$ when $m \geq 6$, we get a contradiction; the cases of $m \leq 5$ will be examined individually.

For an element $g \in G$, let $R(g)$ denote the number of ordered pairs $(a_i,a_j) \in A^2$ for which $g=a_i-a_j$; we say that $R(g)$ is the {\em representation number} of $g$.  Below we prove that $R(g) \in \{0,1,2,m\}$ for each $g \in G$ and, in fact, count the number of elements for each possible representation number.  

Clearly, elements not in $A-A$ have representation number 0, and $R(0)=m$.

Suppose now that $g \in A-A$, $g \neq 0$; we will prove that $R(g)=1$ or $R(g)=2$.  Let's assume indirectly that $R(g) \geq 3$, and let 
$$g=a_1-a_2=a_3-a_4=a_5-a_6$$ for some $a_1, \dots, a_6 \in A$ such that $a_1,a_3,a_5$ are pairwise distinct and $a_2,a_4,a_6$ are pairwise distinct.  (Note that if $a_1=a_3$, for example, then $a_2=a_4$ and thus those two representations of $g$ are identical.)  Our equations imply that $a_1+a_4=a_2+a_3$, $a_1+a_6=a_2+a_5$, and $a_3+a_6=a_4+a_5$.  Since $A$ is a perfect restricted 2-basis, we must have $a_1=a_4$ or $a_2=a_3$; $a_1=a_6$ or $a_2=a_5$; and $a_3=a_6$ or $a_4=a_5$.  In particular, note that if $a_1=a_4$, then we cannot have $a_1=a_6$ or $a_4=a_5$, since that would imply that $a_4=a_6$ or $a_1=a_5$, respectively, which we can't have.  Therefore, there are only two possibilities: $a_1=a_4$, $a_2=a_5$, and $a_3=a_6$; or $a_2=a_3$, $a_1=a_6$, and $a_4=a_5$.  In the first case, we have
$$3g=(a_1-a_2)+(a_3-a_4)+(a_5-a_6)= (a_1-a_2)+(a_3-a_1)+(a_2-a_3)=0,$$ and in the second case 
$$3g=(a_1-a_2)+(a_3-a_4)+(a_5-a_6)= (a_1-a_2)+(a_2-a_4)+(a_4-a_1)=0,$$ and, since $g \ne 0$, this implies that the order of $G$ must be divisible by 3. However, 
$$|G|={m \choose 0} + {m \choose 1} + {m \choose 2} = \frac{m^2+m+2}{2},$$ which is not divisible by 3 for any $m$.  We arrived at a contradiction, and thus every nonzero element of $A-A$ has representation number equal to either 1 or 2.  

We now identify the elements of $A-A$ that have representation number 2.  We call two elements $a_i, a_j \in A$ {\em twins} if $2a_i=2a_j$ or, equivalently, if $a_i - a_j \in L \setminus \{0\}$.  More specifically, we say that $a_i$ and $a_j$ are $g$-twins for some $g \in L \setminus \{0\}$ when $a_i=a_j+g$ (in which case we also have $a_j=a_i+g$).  Clearly, if $a_i$ and $a_j$ are twins, then $R(a_i-a_j)=2$ since $a_i-a_j=a_j-a_i$.  Now we show that there cannot be more than one $g$-twin pair for any $g \in L \setminus \{0\}$; in other words, if $a_i$ and $a_j$ are $g$-twins and $a_i'$ and $a_j'$ are $g$-twins, then $\{a_i,a_j\}=\{a_i',a_j'\}$.  Indeed if $a_i - a_j=a_i'-a_j'=g$, then adding the equations results in $a_i+a_i'=a_j+a_j'$ which, since $A$ is a perfect restricted 2-basis, can only happen if $a_i=a_i'$ (but then $a_j=a_j'$), or $a_j=a_j'$ (but then $a_i=a_i'$), or $a_i=a_j'$ and $a_j=a_i'$.  In any case, $\{a_i,a_j\}=\{a_i',a_j'\}$, as claimed, and this implies that there are at most $|L|-1$ pairs of twins in $A$.  

Next, we examine the {\em awesum} elements of $A$, that is, elements $a_i$ for which $2a_i \in 2 \hat{\;}A$.  If $a_i \in A$ is awesum, then there exists an unordered pair of distinct elements $a_{i'}, a_{i''} \in A$ with $2a_i=a_{i'} + a_{i''};$ we call $a_{i'}$ and $ a_{i''}$ the {\em summands} of $a_i$.  Note that the summands of an element $a_i$ cannot equal $a_i$ since that would imply that the two summands are equal.  Now if $a_i$ is an awesum element of $A$ with summands $a_{i'}$ and $ a_{i''}$, then $R(a_i-a_i')=R(a_i-a_i'')=2$ since $a_i-a_i'=a_i''-a_i$ and $a_i-a_i''=a_i'-a_i$.  

We are now ready to calculate $|A-A|$ as the number of ordered pairs $(a_i,a_j)$ that yield distinct elements of $A-A$.  Since $R(0)=m$ and $R(g) \in \{1,2\}$ for each nonzero element of $A-A$, we have $$|A-A|=m^2-(m-1)-\alpha,$$ where $\alpha$ is the number of elements of $A-A$ whose representation number equals 2.  Suppose that $g \in A-A$ has $R(g)=2$ and that $g=a_1-a_2=a_3-a_4$ for some $a_1,a_2,a_3,a_4 \in G$ with $a_1 \neq a_3$ and $a_2 \neq a_4$.  This implies that $a_1+a_4=a_2+a_3$, so we have three possibilities: (i) $a_1=a_4$ and $a_2=a_3$, in which case $a_1$ and $a_2$ are twins; (ii) $a_1=a_4$ and $a_2 \neq a_3$, in which case $a_1$ is awesum and $a_2$ is one of its summands; or (iii) $a_1 \neq a_4$ and $a_2=a_3$, in which case $a_2$ is awesum and $a_1$ is one of its summands.  Thus $\alpha$ equals the number of twin pairs in $A$ (denoted by $\beta$) plus two times the number of awesum elements of $A$ (denoted by $\gamma$): 
$$|A-A|=m^2-(m-1)-\beta-2\gamma.$$  

Finally, we count the number of elements in $(A-A) \cup A$ by subtracting $|(A-A) \cap A|$ from $|A-A|+|A|$.  Suppose that $a=a'-a''$ for some $a, a', a'' \in A$ or, equivalently, that $a+a''=a'$.  Since $A$ and $2 \hat{\;}A$ are disjoint, this implies that $a=a''$ and so $2a \in A$.  But then  $2a \not \in 2 \hat{\;}A$, so $a$ is not awesum, and thus we have at most $m-\gamma$ elements of $A$ that are in $A-A$.  Putting it all together, we get
\begin{eqnarray*} |(A-A) \cup A| & = & |A-A|+|A|-|(A-A) \cap A| \\
& \geq & m^2-(m-1)-\beta-2\gamma + m-(m-\gamma) \\
& = &  m^2-m+1-\beta-\gamma;
\end{eqnarray*} 
and since $\beta \leq |L|-1$ and $\gamma \leq m$, we arrive at $$|(A-A) \cup A| \geq m^2-2m+2-|L|,$$ as claimed.  

Now since $(A-A) \cup A$ is a subset of $G$, this implies that the order $n$ of $G$ satisfies
$$n=\frac{m^2+m+2}{2} \geq m^2-2m+2-|L|.$$  Recall that $L$ is a subgroup of $G$ and its order is a power of 2.  As here we treat the case when $G$ is not isomorphic to an elementary abelian 2-group, we must have $|L| \leq n/2$, which then implies that $m \leq 10$.  As we mentioned above, $n$ is not divisible by 3 for any value of $m$, and we also find that it is never divisible by 5 either, hence $|L| \neq n/3, n/5,$ or $n/6$.  Furthermore, the values of $n= (m^2+m+2)/2$ are not powers of 2 for any $6 \leq m \leq 10$.  Therefore, we have $|L| \leq n/7$, from which we get $m \leq 5$.

The cases of $m \leq 5$ are easy to handle: since $G$ is not isomorphic to an elementary abelian 2-group, we only need to consider $G \cong \mathbb{Z}_4$ for $m=2$; $G \cong \mathbb{Z}_7$ for $m=3$; $G \cong \mathbb{Z}_{11}$ for $m=4$; and $G \cong \mathbb{Z}_{16}, \; \mathbb{Z}_{2} \times \mathbb{Z}_8, \; \mathbb{Z}_{4}^2,$ or $ \mathbb{Z}_{2}^2 \times \mathbb{Z}_4$ for $m=5$.  One can easily verify that:
\begin{itemize}
  \item in $\mathbb{Z}_4$, the set $\{1,2\}$ is a perfect restricted 2-basis;
  \item in $G=\mathbb{Z}_7$, the set $\{1,2,4\}$ is a perfect restricted 2-basis; and
  \item in $G=\mathbb{Z}_{2}^2 \times \mathbb{Z}_4$, the set $\{(0, 0, 2), (0, 0, 3), (0, 1, 1), (1, 0, 1), (1, 1, 1)\}$ is a perfect restricted 2-basis.
\end{itemize}
(In Section 4 we will see how each of these three constructions generalizes for arbitrary $s$.)    However, we find -- using the computer program \cite{Ili:2017a} -- that $\mathbb{Z}_{11}, \; \mathbb{Z}_{16}, \; \mathbb{Z}_{2} \times \mathbb{Z}_8$ and $\mathbb{Z}_{4}^2$ have no perfect restricted 2-bases.  This completes our proof of Theorem \ref{CoHo-3}.

\section{Further results and open questions}

In this section we generalize the perfect restricted $2$-bases we exhibited above for arbitrary values of $s$ and present some open questions.  

Suppose that a finite abelian group $G$ has a perfect restricted $s$-basis of size $m$.  When $m \le s$, then for the order of $G$ we have
$$|G|=\sum_{h=0}^s {m \choose h} = 2^m.$$

Let $G$ be an abelian group of order $2^m$; we then may assume that there exist positive integers $k_1, \dots, k_r$ for which $$G = \mathbb{Z}_{2^{k_1}} \times \mathbb{Z}_{2^{k_2}} \times \cdots \times \mathbb{Z}_{2^{k_r}}$$ with $k_1 + \cdots + k_r=m$.
For each $i=1,\dots,r$, consider the $k_i$ elements of $G$ whose components are all 0, except for their $i$-th components, which are equal to $1,2, \dots, 2^{k_i-1}$, respectively.  Then it is easy to see that $m$-subset of $G$ consisting of these elements is a perfect restricted  $m$-basis (and thus a perfect restricted  $s$-basis for any $s \geq m$) in $G$.  Therefore, the case of $m \leq s$ is settled:

\begin{prop}

Suppose that $m$ and $s$ are positive integers with $m \leq s$.  Then $G$ has a perfect restricted $s$-basis of size $m$ if, and only if, $G$ has order $2^m$.

\end{prop}

Next, we consider $m=s+1$, in which case $$|G|=\sum_{h=0}^s {m \choose h} = 2^{s+1}-1.$$  As it is easy to see that the set $\{1, 2, , \dots, 2^s\}$ is a  
perfect restricted $s$-basis in the cyclic group of order $2^{s+1}-1$, we have:

\begin{prop}

For every positive integer $s$, the cyclic group of order $2^{s+1}-1$ has a perfect restricted $s$-basis of size $s+1$.

\end{prop}

The question then becomes:

\begin{que}
Are there any noncyclic groups that have a perfect restricted $s$-basis of size $s+1$?

\end{que}

As we have not found any perfect restricted $s$-bases of size more than $s+1$ until reaching size $2s+1$, we ask:

\begin{que}
Are there any groups that have a perfect restricted $s$-basis of size $m$ with $s+2 \leq m \leq 2s$?

\end{que}

This brings us to the case of $m=2s+1$, for which we have
$$\sum_{h=0}^s {m \choose h} = \frac{1}{2} \cdot \sum_{h=0}^{2s+1} {2s+1 \choose h} = 2^{2s}.$$

We found two different groups with perfect restricted $s$-bases of size $2s+1$.  One of them is the elementary abelian 2-group of rank $2s$.  For $i=1,\dots, 2s$, let ${\bf e}_i$ denote the $i$-th unit vector of $\mathbb{Z}_2^{2s}$, that is, the element with a 1 in the $i$-th component and 0 everywhere else, and let ${\bf e}$ be their sum.  Then the $2s$ unit vectors generate each element of the group with at most $s$ nonzero components, and the unit vectors added to ${\bf e}$ generate the rest; each element is generated only once.  We thus have:

\begin{prop}

For every positive integer $s$, the elementary abelian $2$-group of rank $2s$ has a perfect restricted $s$-basis of size $2s+1$.

\end{prop}

The other group of order $2^{2s}$ with a perfect restricted $s$-basis that we found is $\mathbb{Z}_2^{2s-2} \times \mathbb{Z}_4$.  Similarly as above, for $i=1,\dots, 2s-2$ we let ${\bf e}_i$ denote the element of $G$ with a 1 in the $i$-th component and 0 everywhere else, and we let ${\bf e}$ be their sum.  Furthermore, we set ${\bf f}$ equal to the vector with 0s in its first $2s-2$ components and 1 in its last.  We now verify that the set  
$$A=\{{\bf e}_i + {\bf f} \mid i =1,\dots,2s-2\} \cup \{2{\bf f}, 3 {\bf f}, {\bf e}+{\bf f}\}$$ is a perfect restricted $s$-basis in $G=\mathbb{Z}_2^{2s-2} \times \mathbb{Z}_4$.  Since $|A|=2s+1$ and $|G|=2^{2s}$, it suffices to prove that all elements of $G$ can be written as the sum of at most $s$ distinct terms from $A$. 

For an integer $0 \le k \le 2s-2$, we let $S_k$ denote the collection of group elements that have exactly $k$ of their first $2s-2$ components equal to 1, and we write $$B=\{{\bf e}_i + {\bf f} \mid i =1,\dots,2s-2\}$$ for short.  
Now let $g \in S_k$ be an arbitrary element.  Observe that 
$$S_k=\left(k \hat{\;} B \right) \cup \left(k \hat{\;} B + 2{\bf f}\right)  \cup \left(k \hat{\;} B + 3{\bf f} \right) \cup \left(k \hat{\;} B + 2{\bf f} + 3{\bf f} \right) , $$  
and thus if $k \leq s-2$, then $g \in \cup_{h=0}^s h \hat{\;} A.$  
Replacing $k \hat{\;} B $ above by $(2s-2-k) \hat{\;} B +({\bf e}+{\bf f})$ allows us to conclude that $g \in \cup_{h=0}^s h \hat{\;} A$ whenever $2s-2-k$ is at most $s-3$, that is, when $k \geq s+1$.

This leaves us with the cases of $k=s-1$ and $k=s$, for which we see that
$$S_{s-1}=\left((s-1) \hat{\;} B \right) \cup \left((s-1) \hat{\;} B + 2{\bf f}\right)  \cup \left((s-1) \hat{\;} B + 3{\bf f} \right) \cup \left((s-1) \hat{\;} B + ({\bf e} + {\bf f}) \right)$$ and
$$S_{s}=\left(s \hat{\;} B \right) \cup \left((s-2) \hat{\;} B + ({\bf e} + {\bf f})\right)   \cup \left((s-2) \hat{\;} B + ({\bf e} + {\bf f})+ 2{\bf f}\right)  \cup \left((s-2) \hat{\;} B + ({\bf e} + {\bf f})+ 3{\bf f}\right),$$ with which all cases are covered.

Therefore, we have the following result:

\begin{prop}

For every positive integer $s$, the group $\mathbb{Z}_2^{2s-2} \times \mathbb{Z}_4$ has a perfect restricted $s$-basis of size $2s+1$.

\end{prop}  

\begin{que}
Are there other groups besides those isomorphic to  $\mathbb{Z}_2^{2s} $ or $\mathbb{Z}_2^{2s-2} \times \mathbb{Z}_4$ that are of order $2^{2s}$ and have a perfect restricted $s$-basis of size $2s+1$?

\end{que}

We also could not find any cases with $m$ greater than $2s+1$, so ask:

\begin{que}
Are there any groups that have a perfect restricted $s$-basis of size $m$ with $m \geq 2s+2$?

\end{que}

We close by mentioning that while Theorems \ref{CoHo-1} and \ref{CoHo-2} are easily seen to be equivalent, a similar statement does not appear to hold for restricted sums.  In particular, we do not see a correspondence between perfect restricted $s$-bases, considered above, and subsets that generate every group element uniquely as the sum of {\em exactly} $s$ distinct elements.  This then raises the following general question:

\begin{que}
Given a positive integer $s$, characterize all finite abelian groups $G$ that have a subset $A$ for which each element of $G$ arises as a unique sum of exactly $s$ distinct elements of $A$.

\end{que}

\end{document}